\input amstex
\documentstyle{amsppt}

\topmatter

\title
Variations on a Conjecture of Halperin
\endtitle

\author
Gregory Lupton
\endauthor

\abstract
Halperin has conjectured that the Serre spectral
sequence of any fibration that has fibre space a certain kind
of elliptic space should collapse at the $E_2$-term.
In this paper we obtain an equivalent phrasing of this
conjecture, in terms of formality relations between base
and total spaces in such a fibration (Theorem 3.4).  Also, we obtain
results on relations between various numerical
invariants of the base, total and fibre spaces in
these fibrations. Some of our results give weak versions of
Halperin's conjecture (Remark 4.4 and Corollary 4.5).
We go on to establish some of these
weakened forms of the conjecture (Theorem 4.7).
In the last section, we discuss
extensions of our results and suggest some possibilities for future
work.
\endabstract

\address
Department of Mathematics,
Cleveland State University,
2400 Euclid Avenue,
Cleveland OH 44115 U.S.A.
\endaddress

\email
Lupton\@math.csuohio.edu
\endemail

\thanks
This paper was written whilst the author was a guest at
the Max-Planck-Institut f\"ur Mathematik.
The work was begun whilst the author was a
visitor at the Universit\'e Catholique de
Louvain.  Thanks to both institutions for
their support.
\endthanks

\endtopmatter

\document

\heading \S 1 --- Introduction
\endheading

We begin with a description of the conjecture referred to in the
title.  In this paper, all spaces are simply connected CW
complexes and are of finite type over ${\Bbb Q}$, i.e., have
finite-dimensional rational homology groups.  A fibration $F
\overset j \to \longrightarrow E \overset p \to \longrightarrow B$
is said to be {\it totally non-cohomologous to zero\/}
 (abbreviated TNCZ) if
the induced homomorphism $j^*: H^*(E; {\Bbb Q}) \to H^*(F; {\Bbb
Q})$ is onto. This is a very strong condition to place on a
fibration. It is equivalent to requiring that the Serre spectral
sequence (for cohomology with rational coefficients) collapse at
the $E_2$-term (cf. \cite{McC,Th.5.9}). In this case there is an
isomorphism $H^*(E; {\Bbb Q}) \cong H^*(B; {\Bbb Q})\otimes H^*(F;
{\Bbb Q})$ of $H^*(B; {\Bbb Q})$-modules.  Thus a TNCZ fibration
is somewhere between being trivial from the rational homology
point of view and being trivial from the rational cohomology
algebra point of view (cf. Example 1.2).

In the sequel we focus on certain fibre spaces $F$ that satisfy
the following conditions:
\roster
\item $H^*(F; {\Bbb Q})$ is finite-dimensional.
\item $\pi_*(F)\otimes{\Bbb Q}$ is finite-dimensional.
\item The Euler characteristic of $F$, i.e.,
$\sum_i (-1)^i {\text{dim}}\big(H^i(F; {\Bbb Q})\big)$, is
positive.
\endroster
A space that satisfies (1) and (2) is called (rationally) {\it
elliptic\/}.  See \cite{Ha$_1$}, \cite{F\'e,Ch.5} or \cite{Au} for
a discussion of these spaces.  It is known that elliptic spaces
have non-negative Euler characteristic \cite{Ha$_1$}.  So
condition (3) further restricts $F$ to being one of two types of
elliptic space.  We often refer to a space that satisfies
conditions (1)--(3) as a {\it positively elliptic\/} space.
However, we also refer to such spaces in the long-hand, as
`elliptic with positive Euler characteristic', particularly when
stating results.

The conjecture of the title, with which we are concerned, is as
follows:
\proclaim{1.1 Conjecture (Halperin)}  Let $F$ be elliptic with
positive Euler characteristic. Then any fibration $F \to E \to B$
is TNCZ.
\endproclaim
This conjecture has been established in various cases, but in
general it remains open.  Some results that concern it are
mentioned later in the introduction.

We point out that this paper does not resolve Conjecture 1.1, not
even in special cases!  Rather, as the title suggests, we are
concerned here with variations on the theme provided by the
conjecture. These variations come about by considering
consequences of Conjecture 1.1, assuming it to be true.  The
motivation is two-fold:  First, it is hoped to open up new lines
of approach to the conjecture itself.  Second, by looking at such
consequences one can obtain weak versions of the conjecture.  On
the one hand, these weak versions might prove more tractable than
the original.  On the other hand, they should lead to a fuller
understanding of the conjecture.

Although we consider consequences of Conjecture 1.1,
 some of our results
are independent of the status of this conjecture and furthermore
are interesting in their own right.  For instance, Theorem 4.7
specializes to obtain the following result: If $F \to E \to
S^{2n+1}$ is a fibration with fibre a positively elliptic space
and base an odd sphere, then ${\text{cat}}_0(E) =
{\text{cat}}_0(F) + 1$.  This result establishes a weak form of
Conjecture 1.1.  But also, for instance, it can be viewed as a
strong form of Ganea's conjecture in the (very) restricted
circumstances to which it applies.

Next, we outline the contents of the paper. This introductory
section continues with a discussion of positively elliptic spaces
and some of their properties. We go on to discuss models of
rational fibrations, the main technical tool that we use.  The
introduction finishes with a brief summary of some results on
Conjecture 1.1 and some notational conventions. Section 2 is a
short technical section, although in Theorem 2.2 we obtain a very
strong consequence of Conjecture 1.1. In Section 3 we relate the
formality of $E$ and $B$, for a class of fibrations $F \to E \to
B$ including those to which Conjecture 1.1 applies. In Proposition
3.2, for example, we show $B$ formal implies $E$ formal, under the
hypothesis that Conjecture 1.1 is true. We also obtain an
equivalent phrasing of Conjecture 1.1, in Theorem 3.4. In Section
4 we consider some numerical rational homotopy invariants.  Under
the hypothesis that Conjecture 1.1 is true, we obtain inequalities
that relate the values of these invariants on base, total and
fibre spaces of a suitable fibration (Remark 4.4).  These
inequalities can therefore be viewed as weak versions of
Halperin's conjecture.  We go on to establish these weakened forms
of the conjecture in certain restricted circumstances:  For some
of the invariants, we obtain complete results in case the base
space is a wedge of odd-dimensional spheres (Theorem 4.7). In the
last section, we discuss extensions of our earlier results. Here
we suggest various directions for future work, in part by offering
specific questions on these topics.

The spaces $F$ that feature in Conjecture 1.1 are clearly of a
very restricted kind. We continue with a discussion of some of
their properties. We have characterized them by conditions
(1)--(3) above. However, results of Halperin \cite{Ha$_1$} allow
for various characterizations.  Halperin shows that, for an
elliptic space $F$, the three conditions of positive Euler
characteristic, $\chi_{\pi}(F) = 0$ and $H^{\text{odd}}(F;{\Bbb
Q}) = 0$ are equivalent.  Here, $\chi_{\pi}(F)$ denotes the
so-called {\it homotopy Euler characteristic\/} of $F$.  This is a
number defined for any space that has finite-dimensional rational
homotopy by $\chi_{\pi}(F) := \sum_i (-1)^i
{\text{dim}}\big(\pi_i(F)\otimes{\Bbb Q}\big)$. The cohomology
algebra of a positively elliptic space $F$ is zero in odd degrees
and has a presentation of the form
$$H^*(F; {\Bbb Q}) \cong { {\Bbb Q}[x_1,\ldots,x_n] \over
(R_1,\ldots,R_n)}$$
with relations generated by a maximal regular sequence $\{
R_1,\ldots,R_n \}$ in the polynomial algebra ${\Bbb
Q}[x_1,\ldots,x_n]$. Here, the relations $R_1, \ldots, R_n$ need
not be homogeneous (length) polynomials. The {\it minimal model\/}
of such a space is of a particularly restricted form. Recall that
the minimal model of a space $X$ is a differential graded
(henceforth DG) algebra ${\Cal M}_X, d_X$, that as a graded
algebra is a free graded commutative algebra (polynomial on even
degree generators and exterior on odd degree generators).  Also,
its (degree $+1$) differential $d_X$ is decomposable, in the sense
that it induces the trivial differential after passing to the
quotient module of indecomposables, i.e., it has zero linear part.
See \cite{Gr-Mo}, \cite{Ha$_3$} and \cite{Ta} for the basics of
minimal models and their use in rational homotopy theory.  The
book by F\'elix \cite{F\'e} contains more recent material and
references. Condition (1) above implies that the minimal model has
finite-dimensional cohomology, as the cohomology of the minimal
model is identified with that of the space. Condition (2)
translates into the condition that the minimal model be
finitely-generated as a free graded algebra, since the algebra
generators of the minimal model are identified, as a graded vector
space, with the rational homotopy groups of the space.  For an
elliptic space, condition (3) greatly restricts the form of the
minimal model further. It implies, for instance, that up to
isomorphism the model is a {\it pure\/} model \cite{Ha$_1$}. This
is to say that it has the form
$${\Cal M}_F, d_F =
\Lambda(V^{\text{even}})\otimes\Lambda(V^{\text{odd}}), d_F$$
with $d_F(V^{\text{even}}) = 0$ and $d_F(V^{\text{odd}}) \subseteq
\Lambda(V^{\text{even}})$. For an elliptic space, condition (3)
further implies that the minimal model has the same number of even
degree generators as odd degree generators. In symbols, this means
$\text{dim}(V^{\text{even}}) = \text{dim}(V^{\text{odd}})$ and
this fact corresponds to the condition that $\chi_{\pi}(F)=0$.  We
state one more property of these remarkable spaces.  We have said
that the cohomology algebra of a positively elliptic space is zero
in odd degrees.  In fact, more is true.  Suppose $\Lambda V, d$ is
any pure model, as above.  We place a second grading on $\Lambda
V$ by setting $\big(\Lambda V\big)_k = \Lambda V^{\text{even}}
\otimes \Lambda^k V^{\text{odd}}$ for $k \geq 0$.  Since
$d(V^{\text{even}}) = 0$ and $d(V^{\text{odd}}) \subseteq
\Lambda(V^{\text{even}})$, the differential $d$ decreases second
degree by 1 and so $\Lambda V, d$ is a bigraded DG algebra.  This
second grading passes to cohomology.  Now, if $\Lambda V$ is the
model of a positively elliptic space, then we have $H_{+}(\Lambda
V, d) = 0$ (see \cite{Ha$_1$,Th.2}).

It is worthwhile remarking that, despite the highly restricted
form of $F$, there are many examples of such spaces, some of which
correspond to spaces familiar to topologists and geometers:
Even-dimensional spheres, complex projective spaces, Grassmann
manifolds and in general homogeneous spaces $G/H$ with $G$ a
compact, connected Lie group and $H$ a closed subgroup of maximal
rank are all examples of positively elliptic spaces. Further,
given any algebra presented as above, there is some space $F$,
necessarily positively elliptic, that realizes the algebra as its
rational cohomology algebra.

Next, we survey some material on rational fibrations and their
minimal models.  For a fuller discussion see \cite{Ha$_2$} or
\cite{F\'e}. Consider a sequence of DG algebra maps of the form
$$B, d_B  @>i>>  B\otimes\Lambda V, D  @>q>>  \Lambda V, d $$
in which $i \: B, d_B \to B\otimes\Lambda V, D$ is the inclusion
and $q \: B\otimes\Lambda V, D \to \Lambda V, d$ is the projection
onto the quotient DG algebra of $B\otimes\Lambda V$ by the ideal
generated by $B^+$. This sequence is called a {\it KS-extension\/}
of $B, d_B$ (for Koszul-Sullivan) if there is a well-ordered basis
$\{v_{\alpha}\}_{\alpha \in I}$ of $V$ such that, for each $\alpha
\in I$, $D(1\otimes v_{\alpha}) \in B \otimes
\Lambda(V_{<\alpha})$. Here $V_{<\alpha}$ denotes the subspace of
$V$ generated by basis elements $\{v_{\beta} \mid \beta <
\alpha\}$.  Such a well-ordered basis is referred to as a {\it
KS-basis\/} for the extension. In all cases of interest to this
paper, the quotient $\Lambda V, d$ is minimal in the sense
described earlier, i.e., it is free with decomposable
differential.

Next consider a map $p: E \to B$ of 1-connected spaces.  There is
a corresponding DG algebra map of minimal models,  ${\Cal M}(p) \:
{\Cal M}_B, d_B \to {\Cal M}_E, d_E$.  Now any map of DG algebras
can be converted into a KS-extension \cite{F\'e,p.26}. This
process is analogous to that of writing a map of spaces as a
fibration up to equivalence, and results in the following diagram:
$$ \CD {\Cal M}_B, d_B @>{{\Cal M}(p)}>> {\Cal M}_E, d_E  \\ @|
@V{\simeq}V{\phi}V  \\ {\Cal M}_B, d_B @>i>>  {\Cal
M}_B\otimes\Lambda V, D @>q>> \Lambda V, d \\
\endCD
$$
Here, the DG algebra map $\phi$ is a {\it quasi-isomorphism\/},
i.e., it induces an isomorphism on cohomology. A fundamental
result of rational homotopy theory asserts that if $p: E \to B$ is
a fibration of 1-connected spaces with fibre $F$, then $\Lambda V,
d$ is a minimal model of the fibre (cf.{}  \cite{Ha$_2$,Th.4.6} or
\cite{F\'e,Th.2.3.3}).  More generally, we call a sequence of
1-connected spaces $F \to E \to B$ a {\it rational fibration\/}
if, after forming the KS-model ${\Cal M}_B, d_{B} \to {\Cal
M}_B\otimes\Lambda V, D \to \Lambda V, d$ of $E \to B$, the
quotient $\Lambda V, d$ is a minimal model of $F$. In this case,
we refer to the KS-extension ${\Cal M}_B, d_{B} \to {\Cal
M}_B\otimes\Lambda V, D \to \Lambda V, d$ as the minimal model of
the rational fibration. Note that the differential $D$ may have a
non-trivial linear part, even though we use the terminology
`minimal' for such a model (cf. Example 1.2 below). Now suppose $F
\to E \to B$ is a rational fibration with minimal model ${\Cal
M}_B, d_{B} \to {\Cal M}_B\otimes\Lambda V, D \to \Lambda V, d$.
Then the `fibre inclusion' $j: F \to E$ is modeled by the
projection $q\: {\Cal M}_B\otimes\Lambda V, D \to \Lambda V, d$
and the `fibration' $p: E \to B$ by the inclusion $i \: {\Cal
M}_B, d_B \to {\Cal M}_B\otimes\Lambda V, D$.  In particular, it
follows from the fundamental result about minimal models that
$j^*$ is surjective if and only if $q^*$ is surjective.

It is convenient to allow for some flexibility in modeling a given
fibration.  We say two {\it KS-extensions are quasi-isomorphic \/}
if there is a commutative diagram
$$ \CD B, d_B @>i>> B\otimes\Lambda V,D @>q>> \Lambda V,d \\
@V\phi_1V{\simeq}V   @V\phi_2V{\simeq}V @V\phi_3V{\simeq}V \\ B',
d_{B'} @>i'>>  B'\otimes\Lambda V', D' @>q'>> \Lambda V', d' \\
\endCD
$$
in which  each $\phi_i$ is a quasi-isomorphism.  Notice that in
this case, $q'$ is surjective if and only if $q$ is surjective, so
either KS-extension would be sufficient for determining whether
$j^*$ is surjective. By a {\it model\/} for a rational fibration,
we mean any KS-extension that is quasi-isomorphic to the minimal
model.

There is a standard technique for changing the model of a
fibration, known as {\it change of KS-basis\/}.  We mention it
here for the sake of reference. Suppose that $B, d_B \to
B\otimes\Lambda V, D \to \Lambda V, d$ is a KS-extension with
KS-basis $\{ v_{\alpha}\}$.  Suppose given elements $\eta_{\alpha}
\in B^{+}\otimes\Lambda (V_{<\alpha})$, for each $\alpha \in I$.
Define a map of algebras $\phi: B\otimes\Lambda V \to
B\otimes\Lambda V$ by setting $\phi = \iota$ on $B$, and
$\phi(v_\alpha) = v_{\alpha} + \eta_{\alpha}$ on basis elements of
$V$, then extending to an algebra map. Finally, define a new
differential $D'$ on $B\otimes\Lambda V$ by $D' = \phi^{-1} D
\phi$.  Then we have an isomorphism of KS-extensions
$$ \CD B, d_B @>i>> B\otimes\Lambda V,D @>q>> \Lambda V,d \\ @|
@A{\phi}AA  @| \\ B, d_{B} @>i>>  B\otimes\Lambda V, D' @>q>>
\Lambda V, d \\
\endCD
$$
In practice, the new differential $D'$ will have some simpler
form, or display some desired property.  For instance, if we
change basis in such a way that $v_{\alpha} + \eta_{\alpha}$ is a
$D$-cocycle for some $\alpha$, then $D'(v_{\alpha}) = 0$. Another
way of obtaining different models for a fibration is to use the
pushout in the context of KS-extensions. Suppose we have a
KS-extension $B, d_B \to B\otimes\Lambda V, D \to \Lambda V, d$
and a map of DG algebras $\phi \: B, d_B \to B', d_{B'}$. We form
the pushout as described in \cite{Ba,p.66}, for example. (In
Baues' terminology, the inclusion $i$ is a cofibration.) This
gives the following pushout diagram:
$$ \CD B, d_B  @>i>> B\otimes\Lambda V, D \\ @V{\phi}VV
@V{\bar\phi}VV  \\ B', d_{B'} @>{\bar i}>> B'\otimes\Lambda V, D'
\\
\endCD
$$
If $\phi$ is a quasi-isomorphism, then so is $\bar\phi$. Also, the
new differential $D'$ projects to the original $d$ on $\Lambda V$.
Hence, if $\phi$ is a quasi-isomorphism, we obtain a
quasi-isomorphism of KS-extensions as follows:
$$ \CD B, d_B @>i>> B\otimes\Lambda V,D @>q>> \Lambda V,d \\
@V{\phi}V{\simeq}V   @V{\bar\phi}V{\simeq}V @| \\ B', d_{B'}
@>{\bar i}>>  B'\otimes\Lambda V, \delta @>{q'}>> \Lambda V, d \\
\endCD
$$

We illustrate the foregoing discussion of models with an example:

\smallskip

\noindent{\bf 1.2 Example}  There is a fibration $S^2 \to {\Bbb
C}P^3 \to S^4$ (obtained from the Hopf fibration) which has
minimal model
$$\Lambda(w_4, w_7),  d_B \to \Lambda(w_4, w_7)\otimes\Lambda(v_2,
v_3), D \to \Lambda(v_2, v_3) , d.$$
Here, subscripts on generators indicate their degree and the
differentials are given by $d_B(w_4) = 0$, $d_B(w_7) = w_4^2$,
$D(v_2) = 0$ and $D(v_3) = v_2^2 - w_4$.  We obtain $d$ by
projecting onto $\Lambda(v_2, v_3)$.  Notice that $D$ has a
non-trivial linear part. This simple example turns out to have
several interesting features that illustrate quite well the topics
of this paper. In particular, we observe that this fibration is
TNCZ but the cohomology of ${\Bbb C}P^3$ is not isomorphic as an
algebra to the tensor product of the cohomology algebras of $S^4$
and $S^2$.

\smallskip

The model of a rational fibration is the main technical tool that
is used in this paper.  Indeed, the results boil down to algebraic
results about KS-extensions,  which are proven by direct analysis
of the model. When we say `fibration' in the sequel, therefore, we
mean rational fibration, as it is to this class of maps that our
methods apply.

We now mention two ways to re-phrase Conjecture 1.1 that are used
in the sequel. If $F$ is positively elliptic, then its minimal
model is pure.  We extend the notion of pureness to fibrations
with fibre $F$ as follows:

\proclaim{1.3 Definition}  Let $F \to E \to B$ be a fibration in
which $F$ is elliptic with positive Euler characteristic.  The
{\it fibration is pure\/} (as a fibration) if it has minimal model
${\Cal M}_B, d_B \to {\Cal M}_B\otimes\Lambda V, D \to \Lambda V,
d$ in which $D(V^{\text{even}}) = 0$ and $D(V^{\text{odd}})
\subseteq {\Cal M}_B\otimes \Lambda(V^{\text{even}})$.
\endproclaim

Our first rephrasing of Conjecture 1.1 is given in the following:

\proclaim{1.4 Theorem \cite{Th$_1$}} For a given fibration $F \to
E \to B$, in which $F$ is elliptic with positive Euler
characteristic, the following are equivalent: \roster
\item The fibration is TNCZ.
\item The fibration is pure.
\endroster
\endproclaim

\noindent{}Halperin's conjecture is therefore equivalent to the
conjecture that each fibration with fibre a positively elliptic
space is a pure fibration. In \cite{Th$_1$}, Thomas uses his
result to show the conjecture is true if the model for $F$ has
${\text{dim}}(V^{\text{even}}) = 1$ or $2$.

Another re-phrasing of Conjecture 1.1 is given by Meier:

\proclaim{1.5 Theorem \cite{Me}}  Let $F$ be elliptic with
positive Euler characteristic. Then the following are equivalent:
\roster
\item Each fibration $F \to E \to B$,
for arbitrary base space $B$, is TNCZ.
\item Each fibration $F \to E \to S^{2n+1}$,
for $n \geq 1$, is TNCZ.
\item The (graded Lie algebra of)
negative-degree derivations of the cohomology algebra of $F$ are
trivial, ${\text{Der}}^{<0}\,H^*(F; {\Bbb Q}) = 0$.
\endroster
\endproclaim

\noindent{}In \cite{Me}, Meier uses his result to establish some
special cases of the conjecture.

In addition to these results, various special cases of the
conjecture have been established. In \cite{Sh-Te}, the conjecture
is shown to be true for $F$ a homogeneous space of the form $G/H$,
where $G$ is a compact, connected Lie group and $H$ is a closed
subgroup of maximal rank. In \cite{Lu},  the result of Thomas is
extended to the case when ${\text{dim}}(V^{\text{even}}) = 3$.  In
\cite{Ma}, it is shown that the class of spaces for which
Halperin's conjecture is true is closed under fibrations.  This
paper of Markl also contains an introduction to Halperin's
conjecture and some other interesting results. The conjecture has
an interpretation in terms of the rational homotopy theory of
${\text{aut}}_1(F)$, the identity component of the monoid of
self-equivalences of $F$.  This is an interesting aspect of the
conjecture with which, however, we are not directly concerned
 here. See the articles \cite{F\'e-Th} and \cite{Me} for
information and references about this.

We end this introductory section with some notation and
conventions.  In general, we adopt the notation of \cite{F\'e} or
\cite{Ha-St}.  We use $V$ or $W$ to denote a positively graded,
rational graded vector space of finite type. We have already used
$\Lambda V$ to denote the free graded commutative algebra on the
graded vector space $V$.  We write $\Lambda^{+}V$ to denote the
elements of positive degree.  In this paper $\Lambda^{+}V$ is the
augmentation ideal of the canonical augmentation of $\Lambda V$.
Also, $\Lambda^{n}V$ denotes the vector space of polynomials of
homogeneous degree $n$, with respect to the grading by word
length.  If $V$ is zero in even degrees, for example, this agrees
with the $n$'th exterior power of $V$.  Then $\Lambda^{\geq n}V$
denotes the ideal in $\Lambda V$ generated by $\Lambda^{n}V$. The
cohomology of a DG algebra $A, \delta$ is denoted $H(A,\delta)$ or
just $H(A)$.  The elements of positive degree are denoted
$H^{+}(A)$.  Other notation and definitions will be given in the
sequel.

\smallskip

\noindent{\bf Acknowledgement:} I thank Octav Cornea, Yves
F\'elix, John Oprea and Jean-Claude Thomas for several very
helpful conversations about this work.

\heading \S 2 --- Two Preliminary Results
\endheading

We present two technical results about models of the fibrations
with which we are concerned. In these results, we focus on the
case in which the base is a wedge of odd-dimensional spheres or a
single such sphere. From Meier's result, cited as Theorem 1.5,
these fibrations are of particular interest from our point of
view.

The following result is not new (cf{.} \cite{Lu}), but is included
for completeness' sake.

\proclaim{2.1 Lemma} Let $F \to E \to S^{2n+1}$ be a fibration in
which $F$ is elliptic with positive Euler characteristic.  Up to
isomorphism, the minimal model $\Lambda(u) \to
\Lambda(u)\otimes\Lambda V, D \to \Lambda V,d$ has $D$
decomposable and $D(V^{\text{even}}) \subseteq
u\cdot\Lambda^{+}(V^{\text{even}})$.
\endproclaim

\demo{Proof} Suppose $\Lambda(u) \to \Lambda(u)\otimes\Lambda V, D
\to \Lambda V,d$ is the minimal model, with $|u| = 2n+1$. The
differential $D$ is decomposable by \cite{Ha$_2$,Th.1.4(iii)}
(cf{.} Prop.4.12 and Rem.4.18 of the same reference).  Now let $v
\in V^{\text{even}}$ and write $D(v) = u\chi_{0} + u\chi_{+}$,
with $\chi_{0} \in (\Lambda^{+}V)_0$ and $\chi_{+} \in
\big(\Lambda V\big)_{+}$.  The subscripts refer to the second
grading of $\Lambda V$ mentioned in the introduction. Applying $D$
again, we obtain $0 = D^2(v) = -ud(\chi_0) - ud(\chi_{+}) =
-ud(\chi_{+})$.  Therefore, $\chi_{+}$ is a $d$-cocycle in
$\big(\Lambda V\big)_{+}$. As stated in the introduction, each
cocycle of positive second degree is a boundary, so there is an
element $\eta \in \Lambda V$ with $d\eta = \chi_{+}$. Hence
$D(u\eta) = -u\chi_{+}$ and we have $D(v+ u\eta) = u\chi_0$.

Now use the change of KS-basis argument mentioned in the
introduction, replacing each KS-basis element $v \in
V^{\text{even}}$ with $v + u\eta$. In the isomorphic KS-extension
obtained as a result, we have $D'(v) = u\chi_0$, as is easily
checked. \hfill$\square$
\enddemo

The next result shows that Halperin's conjecture implies the
strongest possible restriction on a fibration with base a wedge of
odd spheres. It uses the fact that a wedge of spheres is a {\it
formal\/} space. Recall that a space is formal if its minimal
model is determined by its cohomology algebra. Specifically, $X$
is formal if there is a quasi-isomorphism
$$\psi: {\Cal M}_X,d_X \to H^*(X; {\Bbb Q}), 0.$$
We mention in passing that there are many interesting examples of
formal spaces including $H$-spaces and co-$H$-spaces, symmetric
spaces and simply connected, compact K\"ahler manifolds (for this
last assertion, see \cite{D-G-M-S}). A product or wedge of formal
spaces is again a formal space. A positively elliptic space is a
formal space. For further discussion of formal spaces, see Section
3.

\proclaim{2.2 Theorem} Let $F \to E \to B$ be a fibration in which
$F$ is elliptic with positive Euler characteristic and $B$ is
rationally a wedge of odd spheres. If the fibration is TNCZ, then
$E \simeq_{\Bbb Q} B \times F$ and the fibration (rationally) is
trivial.
\endproclaim

\demo{Proof}  Suppose $\Lambda W, d_B \to \Lambda W\otimes\Lambda
V, D \to \Lambda V, d$ is the minimal model of the fibration.  By
the result of Thomas, cited as Theorem 1.4, this model can be
assumed pure. Since $B$ is rationally a wedge of spheres, it is
formal and hence there is a quasi-isomorphism $\psi \: \Lambda W,
d_B \to H(B)$.  Now form the pushout and obtain a quasi-isomorphic
KS-extension $H(B) \to H(B)\otimes\Lambda V, D' \to \Lambda V, d$.
Since the original KS-extension was pure, it follows directly from
the pushout construction that this KS-extension is also pure. The
observation that we have such a model also follows from a close
reading of Thomas' proof of Theorem 1.4.

Now consider the pure model $H(B) \to H(B)\otimes\Lambda V, D \to
\Lambda V, d$. Let $v \in V^{\text{odd}}$ and let $\{b_i\}$ be a
basis for $H^+(B)$. Write the differential $D$ as
$$D(v) = dv + \sum_i\,b_i \Omega_i(v). \eqno (*)$$
For parity of degree reasons each $\Omega_i(v)$ is an odd-degree
element of $\Lambda V$. Since this model is pure,  $D(dv) = 0$.
Therefore, applying $D$ to $(*)$ gives $0 = D^2(v) = -\sum_i b_i
d\Omega_i(v)$. For the  last expression, we use the fact that
products in $H^{+}(B)$ are trivial.  Now we have $d\Omega_i(v) =
0$ for each $i$, so each $\Omega_i(v)$ is a $d$-cocycle of odd
degree, hence of positive second degree in the second grading of
$\Lambda V$ mentioned in the introduction. However, $H_{+}(\Lambda
V) = 0$. Therefore, we can choose elements $\eta_i \in \Lambda V$
for which $\Omega_i(v) = d\eta_i$.  This gives $D\big(\sum_i\, b_i
\eta_i\big) = - \sum_i\, b_i \Omega_i(v)$.

Finally, make a change of KS-basis, replacing each KS-basis
element $v \in V^{\text{odd}}$ with $v + \sum_i\, b_i \eta_i(v)$.
This gives a quasi-isomorphic KS-extension $H(B) \to
H(B)\otimes\Lambda V, D' \to \Lambda V, d$ in which $D' = 1\otimes
d$, as is easily checked.  The result follows. \hfill$\square$
\enddemo

\smallskip

\noindent{\bf 2.3 Remark.} Notice that we cannot relax the
hypothesis on the base to allow even-dimensional spheres.  Indeed,
Example 1.2 is a non-trivial fibration with positively elliptic
fibre and base an even-dimensional sphere.   For more comments
along these lines, see Section 5.

\heading \S 3 --- Formality and TNCZ Fibrations
\endheading

Halperin's conjecture asserts that a fibration is close to trivial
if the fibre is positively elliptic. So suppose a fibration has
base a formal space and fibre a (formal) positively elliptic
space. Then the conjecture asserts that the total space is close
to a product of formal spaces and therefore close to being formal.
This point of view turns out to yield an equivalent formulation of
the conjecture.

A formal space has minimal model that is more highly structured
than an arbitrary minimal model.  We recall here some properties
of the {\it bigraded model\/} that are used in the sequel. See
\cite{Ha-St} for a full discussion.  The bigraded model of a
formal space is a minimal model $\Lambda V, d$ for which the
vector space of generators has a second grading, $V =
\oplus_{i\geq0}V_i$.  This gives $\Lambda V$ the structure of a
bigraded algebra.  Furthermore, the bigraded model is a bigraded
DG algebra in the sense that $d \: \Lambda V \to \Lambda V$
decreases second grading by exactly 1.  In particular, $d(V_0) =
0$. The second grading therefore passes to cohomology and, as
further properties of the bigraded model, we have $H(\Lambda V,d)
= H_0(\Lambda V,d)$ and $H_{+}(\Lambda V,d) = 0$.  There is a
quasi-isomorphism
--- indeed, the bigraded model proper --- given by
$$\rho \: \Lambda V,d \to H(\Lambda V,d) = H_0(\Lambda V,d)$$
which maps $V_{+}$ to zero and each $v \in V_0$ to the class it
represents in $H(\Lambda V,d)$. A positively elliptic space is
formal and its bigraded model has a simple form.  We identify $V_0
= V^{\text{even}}$ and $V_1 = V^{\text{odd}}$.  There are no
generators of second degree $\geq 2$.  Notice that this accords
with the fact that $H_{+}(\Lambda V,d) = 0$, as we stated in the
introduction. Following \cite{Ha-St}, we adopt the notation
$V_{(n)} = \oplus_{i=0}^n V_i$ and $(\Lambda V)_{(n)} =
\oplus_{i=0}^n (\Lambda V)_i$.

\proclaim{3.1 Proposition} Let $F \to E \to B$ be a fibration in
which $F$ and $B$ are formal. If the fibration is TNCZ, then there
is a model
$$H(B),0  \to H(B)\otimes\Lambda V, D \to \Lambda V, d,$$
in which $\Lambda V, d$ is the bigraded model of $F$ and
$H(B)\otimes\Lambda V,D$ is filtered in the following sense:  Let
$V = \oplus_{i\geq 0} V_i$ be the second grading of the bigraded
model.  Then $D(V_0) = 0$ and for each $i \geq 1$, $D(V_i)
\subseteq H(B)\otimes(\Lambda V)_{(i-1)}$.
\endproclaim

\demo{Proof}  Since the base and fibre are formal, we can suppose
there is a model $H(B) \to H(B)\otimes\Lambda V, \delta \to
\Lambda V, d$ in which $\Lambda V, d$ is the bigraded model of
$F$. Suppose that the fibration is TNCZ.  We begin by showing that
up to isomorphism, $\delta(V_0) = 0$.

Let $v \in V_0$, so that $v$ represents a class $[v]$ in
$H(\Lambda V)$. By assumption, there is some $\delta$-cocycle
$\chi \in H(B)\otimes\Lambda V$ with $q^*([\chi]) = [v]$. It is
easy to see that we may choose the $\delta$-cocycle to be $\chi =
v + \beta$ with $\beta \in  H^{+}(B)\otimes\Lambda V$. Now use the
change of KS-basis argument, replacing each KS-basis element $v
\in V_0$ by $\chi$.  Observe that in the isomorphic model, $H(B)
\to H(B)\otimes\Lambda V,\delta_0 \to \Lambda V, d$, we have
$\delta_0(V_0) = 0$.

Now suppose inductively that for some $n \geq 0$, we have a model
$H(B) \to H(B)\otimes\Lambda V,\delta_n \to \Lambda V, d$ in which
$\delta_n$ is a filtered differential on $H(B)\otimes\Lambda
V_{(n)}$, i.e., $\delta_n(V_i) \subseteq H(B)\otimes(\Lambda
V)_{(i-1)}$ for $i=0,\dots,n$. Let $v \in V_{n+1}$ and write
$$\delta_n(v) = dv + \xi_{(n)} + \xi_{+}\eqno (*_1)$$
for elements $\xi_{(n)} \in H^{+}(B)\otimes(\Lambda V)_{(n)}$ and
$\xi_{+} \in  H^{+}(B)\otimes(\Lambda V)_{\geq n+1}$.

\noindent{\it Claim.\/}  There is an element $\eta \in
H^{+}(B)\otimes\Lambda V$ for which $\delta_n(v) = dv + \xi_{(n)}
+ \delta_n(\eta)$.

\noindent{\it Proof of Claim.\/} Suppose that $\xi_{+} = \xi^m_{+}
+ \xi^{m+1}_{+} +\cdots +\xi^M_{+}$, for some $2 \leq m \leq M$,
with each $\xi^i_{+} \in  H^i(B)\otimes(\Lambda V)_{\geq n+1}$.
Let  $\{b_j\}$ be a basis for $H^m(B)$ and write the `lowest' term
$\xi^m_{+}$ as $\xi^m_{+} = \sum_j\,b_j\otimes\chi_j$ for suitable
elements $\chi_j \in (\Lambda V)_{\geq n+1}$. Applying $\delta_n$
to $(*_1)$ gives
$$0 = \delta_n^2(v) =  \delta_n(dv + \xi_{(n)}) +
 \delta_n(\xi^m_{+}) +
 \sum_{i=m+1}^M  \delta_n(\xi^i_{+}).\eqno(*_2)$$
The induction hypothesis on $\delta_n$ implies that $\delta_n(dv +
\xi_{(n)}) \in H(B)\otimes(\Lambda V)_{(n-1)}$. The ideal in
$H(B)\otimes\Lambda V$ generated by elements of $H(B)$ of degree
at least $m+1$ is $\delta_n$-stable and contains all terms
contributed by the $\delta_n(\xi^i_{+})$ for $i \geq m+1$.
Therefore, in equation $(*_2)$, terms of the form $\sum_j (-1)^{m}
b_j\otimes d(\chi_j)$ are the only contributions to
$H^m(B)\otimes(\Lambda V)_{\geq n}$.  It follows that each
$d(\chi_j) = 0$.  Since this is the bigraded model, in which
$H_{+}(\Lambda V) = 0$, there are elements $\eta_j \in \Lambda V$
such that $(-1)^m\chi_j = d(\eta_j)$, for each $j$. Thus we have
$$\delta_n\big(\sum_j b_j\otimes\eta_j\big) =  \xi^{m}_{+} +
\zeta^{m+1}_{+} + \cdots + \zeta^{M'}_{+},$$
with each $\zeta^{j}_{+} \in  H^j(B)\otimes(\Lambda V)_{\geq
n+1}$,
  for $m+1 \leq j \leq M'$.
Substituting this into $(*_1)$ above gives
$$\delta_n(v) = dv + \xi_{(n)} + \xi'_{+} + \delta_n\big(\sum_j
b_j\otimes \eta_j\big),$$
where $\xi'_{+} = \xi_{+} - \xi^m_{+} - \sum_{j\geq m+1}
\zeta^j_{+} \in H^{\geq m+1}(B)\otimes(\Lambda V)_{\geq n+1}$. An
induction argument on $m$, repeating this last step as necessary,
shows the claim. {\it End of Proof of Claim.\/}

Hence we can make another change of KS-basis, this time replacing
a KS-basis element $v\in V_{n+1}$ by $v - \eta$. In the isomorphic
model $H(B) \to  H(B)\otimes\Lambda V,\delta_{n+1} \to \Lambda V,
d$ we have $\delta_{n+1}(v) = dv + \xi_{(n)}$, with $\xi_{(n)} \in
H^{+}(B)\otimes(\Lambda V)_{(n)}$.  Therefore, $\delta_{n+1}$ is a
filtered differential on $H(B)\otimes \Lambda V_{(n+1)}$.

To finish, use the inductive step just proved to make a global
change of KS-basis, working inductively over $n$.  This results in
a filtered model as required. \hfill$\square$
\enddemo

We now develop the main result of the section. The next
proposition is a little more general than we need in the sequel.
However, it is interesting in its own right.  In it, we only
assume that the fibre space is formal and elliptic.  Recall
Halperin's result \cite{Ha$_1$}, that any elliptic space has
non-negative Euler characteristic. If a space is positively
elliptic, then it is formal.  However, a space that is formal and
elliptic need not have positive Euler characteristic.  For
instance, a product of odd-dimensional spheres is elliptic and
formal, but has Euler characteristic equal to zero.   There are
strong restrictions on a space, however, that follow from the
hypothesis of formal and elliptic.  We do not dwell on this point
here, but simply point out that any space that is formal and
elliptic has a two-stage bigraded model (cf. \cite{F\'e-Ha$_1$}).
This is the feature of these spaces that we use in the result.

\proclaim{3.2 Proposition} Let $F \to E \to B$ be a fibration in
which $F$ is formal and elliptic and $B$ is formal.  If the
fibration is TNCZ, then $E$ is formal also.
\endproclaim

\demo{Proof} From Proposition 3.1 we have a model of the fibration
$$H(B),0  \to H(B)\otimes\Lambda V, D, \to\Lambda V, d$$
in which $D(V_0) = 0$ and for $i\geq 1$, $D(V_i) \subseteq
H(B)\otimes(\Lambda V)_{(i-1)}$. This means, in particular, that
$D(V_1) \subseteq H(B)\otimes(\Lambda V)_{0}$. Since the bigraded
model for $F$ is two-stage, and thus $V = V_0 \oplus V_1$, the
total space $H(B)\otimes\Lambda V, D$ in this model is actually a
bigraded DG algebra.

We now show that, with respect to this bigrading,
$H_{+}\big(H(B)\otimes\Lambda V\big) = 0$. For let $x \in
(H(B)\otimes\Lambda V)_r$ be a $D$-cocycle, for any $r \geq 1$.
We show that $x$ is exact with an argument similar to that used to
show the claim in Proposition 3.1. Write $x$ as $x = x_m + x_{m+1}
+ \cdots + x_M$, for $2 \leq m \leq M$, where each $x_i \in
H^i(B)\otimes(\Lambda V)r$. Let  $\{ b_{j} \}$ be a basis for
$H^m(B)$. Then we can write the `lowest' term $x_m$ as $x_m =
\sum_{j}  b_{j}\otimes\chi_j$ for appropriate terms $\chi_j \in
(\Lambda V)_r$. Now $(D-d)(\chi_j) \in H^{+}(B)\otimes \Lambda V$
and also the ideal $H^{\geq m+1}(B)\otimes \Lambda V$ is
$D$-stable. Therefore, as $x$ is a $D$-cocycle, we have
$$0 = D(x) \equiv (-1)^m \sum_{j} b_{j}\otimes d(\chi_j),$$
where the congruence is modulo the ideal $H^{\geq m+1}(B)\otimes
\Lambda V$.  Therefore, $d\chi_j = 0$ and, since $\Lambda V, d$ is
the bigraded model and since $r$ is positive, there are elements
$\eta_j$ with $d\eta_j = (-1)^m\chi_j$ for each $j$.  This yields
$D\big(\sum_{j} b_{j}\otimes\eta_j\big) \equiv \sum_{j}
b_{j}\otimes\chi_j$, again modulo $H^{\geq m+1}(B)\otimes \Lambda
V$.  Now we can write
$$x = D\big(\sum_{j} b_{j}\otimes\eta_j\big)
 + x'_{m+1} + \cdots + x'_{M'},$$
with each $x'_i \in H^i(B)\otimes(\Lambda V)r$. An induction
argument repeating this argument as necessary obtains the result
that $x$ is $D$-exact.

We have shown that $H\big(H(B)\otimes\Lambda V,\delta\big) =
H_0\big(H(B)\otimes\Lambda V,\delta\big)$.  But now the projection
$$p: H(B)\otimes\Lambda V,\delta \to {H(B)\otimes\Lambda V,\delta
\over \big( \delta  V_1 \big)} =  H_0\big(H(B)\otimes\Lambda
V,\delta\big)$$
is a quasi-isomorphism.  Since $H(B)\otimes\Lambda V,\delta$ is
quasi-isomorphic to the minimal model for $E$, and also $p$ is a
quasi-isomorphism from it to its cohomology, it follows that $E$
is formal. \hfill$\square$
\enddemo

The next result is something of a converse to Proposition 3.2.

\proclaim{3.3 Proposition} Let $F$ be be elliptic with positive
Euler characteristic and let $F \to E \to S^{2n+1}$ be a
fibration. If $E$ is formal, then the fibration is TNCZ.
\endproclaim

\demo{Proof} From Lemma 2.1 the fibration has a model
$$\Lambda(u) \to \Lambda(u)\otimes\Lambda V, D \to \Lambda V,d,$$
in which $\Lambda(u)\otimes\Lambda V, D$ is actually the minimal
model of $E$, and $D(V_0) \subseteq u\cdot\Lambda^{+} V_0$. We
will show that $E$ formal implies $D(V_0) = 0$.  Recall the
characterization of formality given in \cite{D-G-M-S, Th.4.1}.
This says that there is a vector space decomposition $\langle u
\rangle\oplus V \cong C\oplus N$ with $D(C) = 0$, $D: N \to
\Lambda(u)\otimes\Lambda V$ injective and such that any cocycle in
the ideal $I(N)$ of $\Lambda(u)\otimes\Lambda V$ generated by $N$
is exact. Clearly $u \in C$, and we show that $V_0 \subseteq C$,
for any such decomposition.  For suppose not, so that $D(V_0)
\not= 0$. Choose a KS-basis $V_0 = \langle x_1, \ldots, x_n
\rangle$. There is at least one of the $x_i$ with non-zero
differential, so let $r$ be the largest subscript with $D(x_r)
\not= 0$. Note that $x_r \in N$ in the decomposition. Now $D(V_0)$
is contained in the ideal of $\Lambda(u)\otimes\Lambda V$
generated by elements $\{ u x_1, \ldots, u x_{r-1}\}$, which we
denote  $I_{r-1}$. For parity of degree reasons, $(D-d)(V_1)$ is
contained in the ideal generated by $u\cdot V_1$.  Also, we have
$d(V_1) \subseteq \Lambda^{\geq 2}V_0$. Let $J$ denote the ideal
of $\Lambda(u)\otimes\Lambda V$ generated by $u\cdot
V_1+\Lambda^{\geq 2}V_0$.  Then the image of $D$ is contained in
the ideal $I_{r-1} + J$. Now $D(x_r) = u\chi$ for some $\chi$ and
hence we obtain a $D$-cocycle $ux_r \in I(N)$.  This cannot be
exact, as it is not in the ideal $I_{r-1} + J$. This contradicts
the assumptions on the decomposition $C\oplus N$.  Therefore, we
must have $V_0 \subseteq C$ and so $D(V_0) = 0$.  From this it
follows easily that the fibration is TNCZ, because then $q^*\:
H(\Lambda(u)\otimes\Lambda V) \to H(\Lambda V)$ is surjective onto
the generators of $H(\Lambda V)$. \hfill$\square$
\enddemo

Finally, we collect together the preceding results into the main
result of the section.  As we see, we have obtained an equivalent
formulation of Halperin's conjecture.

\proclaim{3.4 Theorem} Let $F$ be elliptic with positive Euler
characteristic.  Then the following are equivalent: \roster
\item Any fibration with fibre $F$ is TNCZ.
\item For any fibration $F \to E \to B$ in
which $B$ is formal, $E$ is formal also.
\item For any fibration $F \to E \to S^{2n+1}$,
the total space $E$ is formal.
\endroster
\endproclaim

\demo{Proof} The implication (1) $\Rightarrow$ (2) follows from
Proposition 3.2. (2) $\Rightarrow$ (3) is obvious, since spheres
are formal spaces. Assume (3).  Then Proposition 3.3 implies that
each fibration $F \to E \to S^{2n+1}$ is TNCZ.  Hence, from
Meier's result (Theorem 1.5), any fibration $F \to E \to B$ is
TNCZ. \hfill$\square$
\enddemo

\heading \S4 --- Numerical Invariants
\endheading

Here we consider some invariants related to the {\it
Lusternik-Schnirelmann category\/}. Recall that this is a
numerical homotopy invariant of a space, defined as one less than
the smallest number of open sets required to cover the space, when
each is contractible in the space. As is usual in rational
homotopy theory, we have `normalised' so that a sphere has
category equal to 1. This invariant and its approximations have
been much studied both in ordinary and rational homotopy theory.
See \cite{Ja} for a recent survey with many references. Here we
focus on four rational homotopy invariants. We define these
invariants and include some discussion, before proceeding to the
results:

{\it (Rational) Cup-length\/}: This is the nilpotency---as an
algebra---of the rational cohomology algebra of a space $X$.  It
is denoted here cup$_0(X)$.  For example we have cup$_0({\Bbb
C}P^n) = n$ for each $n \geq 1$.

{\it (Rational) Toomer's invariant\/}: As in
\cite{F\'e-Ha$_2$,Rem.9.3} we describe this invariant as follows:
Let $\Lambda V,d$ be the minimal model of $X$.  Consider the
projection
$$p_n \: \Lambda V \to {\Lambda V \over \Lambda^{\geq n+1}V}.$$
We obtain our rational invariant, denoted ${\text{e}}_0(X)$, by
setting ${\text{e}}_0(X) \leq n$ if $(p_n)^*$ is injective.  In
other words, ${\text{e}}_0(X)$ is the largest $n$ for which some
non-zero class in $H(\Lambda V)$ is represented by a cocycle in
$\Lambda^{\geq n}V$.

{\it (Rational) Category\/}: This is the Lusternik-Schnirelmann
category of the rationalization of $X$.  We denote it cat$_0(X)$
and, following \cite{F\'e-Ha$_2$,Th.4.7}, we describe it in terms
of the minimal model of $X$:  Set ${\text{cat}}_0(X) \leq n$ if
the above projection $p_n$ makes $\Lambda V$ into a retract of the
quotient $\Lambda V \slash \Lambda^{\geq n+1}V$.

{\it (Rational) Cone-length\/}: This is the least number of steps
required to build the rational homotopy type of a space $X$ as a
succession of cofibration sequences of rational spaces. It is
denoted cl$_0(X)$. Specifically, set ${\text{cl}}_0(X) = 0$ if
$X\simeq_{\Bbb Q}*$ and ${\text{cl}}_0(X) = 1$ if $X$ has the
rational homotopy type of a wedge of spheres.  In general, set
${\text{cl}}_0(X) \leq n$ if there are spaces $X_1, A_1, \dots,
A_{n-1}$, each of which has the homotopy type of a wedge of
rational spheres, and $n-1$ cofibration sequences $A_i \to X_i \to
X_{i+1}$, for $i=1, \dots, n-1$, such that $X_n \simeq_{\Bbb Q}
X$. In \cite{Co} it is shown that ${\text{cl}}_0(X)$ agrees with
the `homotopical nilpotency' of the minimal model of $X$, i.e.,
the least $n$ for which the minimal model is quasi-isomorphic to a
DG algebra that is of nilpotency $n$ as an algebra.  It is in this
latter guise that we meet this invariant here.

We have been a little careless in phrasing the above definitions
by implicitly assuming these invariants finite.  Although the case
when one or other of these is infinite does not seem so
interesting in our context, it is allowed for, where appropriate,
in the following results.

For these invariants, we always have inequalities as follows:
$${\text{cup}}_0(X) \leq {\text e}_0(X) \leq {\text{cat}}_0(X)
\leq  {\text{cl}}_0(X).$$
In the special case that $X$ is a formal space, all these
invariants agree.  In this case we will denote their common value
${\text{nil}}_0(X)$. This usage accords with Cornea's homotopical
nilpotency in this case.  In particular, if $F$ is a positively
elliptic space, then it is formal and we have ${\text{cup}}_0(F) =
{\text e}_0(F) = {\text{cat}}_0(F) = {\text{cl}}_0(F)$, which we
denote by ${\text{nil}}_0(F)$.

We mention some examples to illustrate these invariants:

\noindent{\bf 4.1 Examples} If $X = S^n$, then $X$ is formal, our
four invariants agree and we have ${\text{nil}}_0(X) = 1$. If $X =
{\Bbb C}P^n$, it is likewise formal and ${\text{nil}}_0(X) = n$.
Next, suppose $X = S^2 \vee S^2 \cup_{\alpha} e^5$, where $\alpha
= [\iota_1, [\iota_1,\iota_2]]$, the triple Whitehead product in
$\pi_4(S^2 \vee S^2)$. Then we have  ${\text{cup}}_0(X) = 1$  but
${\text{e}}_0(X) = {\text{cat}}_0(X) = {\text{cl}}_0(X) =2$. A
well-known example of Lemaire-Sigrist, developed by
F\'elix-Halperin (cf. \cite{F\'e-Ha$_2$}), is $X = ({\Bbb C}P^2
\vee S^2)\cup_{\omega} e^7$ for a certain attaching map $\omega$.
This space satisfies ${\text{cup}}_0(X) = {\text{e}}_0(X) = 2$
whilst ${\text{cat}}_0(X) = {\text{cl}}_0(X) =3$. Furthermore,
this example has ${\text{e}}_0(X^n) = 2n$ and ${\text{cat}}_0(X^n)
= 3n$.  This illustrates that ${\text{e}}_0(X)$ can be smaller
than ${\text{cat}}_0(X)$ by an arbitrarily large amount.

These invariants behave quite well for products of spaces. The
product formula ${\text{cup}}_0(X\times Y) = {\text{cup}}_0(X) +
{\text{cup}}_0(Y)$ is well-known. It is easy to see that ${\text
e}_0$ likewise is additive for products
\cite{F\'e-Ha$_2$,Rem.9.3}. Recently, ${\text{cat}}_0$ has been
shown to be additive for products \cite{F\'e-Ha-Le}, and
${\text{cl}}_0$ to be additive at least for products of rational
Poincar\'e duality spaces. Indeed, it has been shown that
${\text{e}}_0(X) = {\text{cat}}_0(X) = {\text{cl}}_0(X)$ whenever
$X$ is a rational Poincar\'e duality space (see \cite{F\'e-Ha-Le}
and \cite{Co-F\'e-Le}). Now Halperin's conjecture asserts that in
certain fibrations the total space is close to being a product of
the base and fibre spaces.  These remarks combine to suggest there
should be good relations between these invariants for base, total
and fibre spaces in such fibrations.  We shall see that this is
the case.

First we give some results that complement one of Jessup
\cite{Je,Prop.3.6}. Let $F \to E \to B$ be a TNCZ fibration with
$F$ formal. Then Jessup's result gives ${\text{cat}}_0(E) \geq
{\text{cat}}_0(B)+{\text{nil}}_0(F)$. Actually, his result applies
a little more generally and was proved for (what was then thought
to be) a different numerical invariant, ${\text{Mcat}}_0$.  The
conclusion for ${\text{cat}}_0$ follows by a result of Hess
\cite{He}, which identified ${\text{cat}}_0$ with
${\text{Mcat}}_0$.

We specialize this result to the following:

\proclaim{4.2 Proposition}  Let $F \to E \to B$ be a fibration in
which $F$ is elliptic with positive Euler characteristic and $B$
is formal.  If the fibration is TNCZ, then $E$ is also formal and
we have ${\text{nil}}_0(E) \geq {\text{nil}}_0(B) +
{\text{nil}}_0(F)$.
\endproclaim

\demo{Proof} The formality of $E$ follows immediately from
Proposition 3.2. Since $E$ and $B$ are both formal,
${\text{nil}}_0 = {\text{cat}}_0$ for these spaces and the
conclusion follows from the result of Jessup mentioned above.
\hfill$\square$
\enddemo

Next, consider the case in which $B$ is not formal. Notice the
following result does not require that $F$ be positively elliptic.

\proclaim{4.3 Proposition}  Let $F \to E \to B$ be a fibration
with $F$ formal.  If the fibration is TNCZ, then we have the
following inequalities: \roster
\item ${\text{cup}}_0(E) \geq {\text{cup}}_0(B) +
{\text{nil}}_0(F)$.
\item ${\text{e}}_0(E) \geq {\text{e}}_0(B) +
{\text{nil}}_0(F)$.
\endroster
\endproclaim

\demo{Proof} Actually, for (1) the hypothesis of formality on $F$
is redundant.  Recall from the introduction, that if a fibration
$F \to E \to B$ is TNCZ, then $H^*(E;{\Bbb Q})\cong H^*(B;{\Bbb
Q})\otimes H^*(F;{\Bbb Q})$ as $H^*(B;{\Bbb Q})$-modules.  Under
this isomorphism, two elements of $1\otimes H^*(F;{\Bbb Q})$
multiply as $(1\otimes x)(1\otimes y) = 1\otimes xy + \chi$, for
some $\chi$ in the ideal generated by $H^{+}(B)$. It follows that
${\text{cup}}_0(E) \geq {\text{cup}}_0(B) + {\text{cup}}_0(F)$.

For (2), we work with a model of the fibration. Suppose that $F$
is formal and the fibration is TNCZ. Then it has a model $\Lambda
W, d_B\to \Lambda W\otimes\Lambda V, D \to \Lambda V, d$ in which
$D(V_0) = 0$. This follows from the argument in the first two
paragraphs of the proof of Proposition 3.1, replacing $H(B)$ there
by $\Lambda W$. So let $\alpha \in \Lambda^{\geq n} V$ be a
cocycle representative for some non-zero class in $H(\Lambda V)$.
Since $F$ is formal, we can suppose that $\alpha \in \Lambda^{\geq
n} V_0$. In our model we have $D(V_0)=0$, so $[\alpha] \in
H(\Lambda W\otimes\Lambda V)$. Furthermore, $\alpha$ is not
$D$-exact, since it is not $d$-exact. Let $\beta\in \Lambda^{\geq
m} W$ be a cocycle representative for some non-zero class of
$H(\Lambda W)$. Since we have $H(\Lambda W\otimes\Lambda V) \cong
H(\Lambda W)\otimes H(\Lambda V)$ as $H(\Lambda W)$-modules, the
product $[\beta][\alpha] = [\beta\alpha]$ is non-zero in
$H(\Lambda W\otimes\Lambda V)$. Now $\alpha \beta\in \Lambda^{\geq
m+n}(W\oplus V)$ so  ${\text{e}}_0(E) \geq m+n$. \hfill$\square$
\enddemo

Example 1.2 illustrates that inequality, as in Proposition 4.3, is
the best that can be hoped for in general.  But see below for
stronger relations in special cases.

If $B$ is not formal in Proposition 4.3, there is no {\it a
priori\/} reason why the invariants ${\text{cup}}_0$,
${\text{e}}_0$ and ${\text{cat}}_0$ for $B$ or $E$ should agree.
Therefore, it can be thought of as giving three distinct necessary
conditions for Halperin's conjecture to be true.  We summarise
this in the following:

\smallskip

\noindent{\bf 4.4 Remark.} Let $F \to E \to B$ be a fibration in
which $F$ is elliptic with positive Euler characteristic. If $B$
is formal and if the fibration is TNCZ, then $E$ is also formal
and ${\text{nil}}_0(E) \geq {\text{nil}}_0(B) +
{\text{nil}}_0(F)$. For general $B$, if the fibration is TNCZ then
we have three inequalities
$$\align {\text{cup}}_0(E) &\geq {\text{cup}}_0(B) +
{\text{nil}}_0(F)\\ {\text{e}}_0(E) &\geq {\text{e}}_0(B) +
{\text{nil}}_0(F) \\ {\text{cat}}_0(E) &\geq {\text{cat}}_0(B) +
{\text{nil}}_0(F).\\
\endalign
$$
Each of these inequalities gives a necessary condition for
Halperin's conjecture to be true.

\smallskip

A much stronger consequence follows if we restrict the base as
follows:

\proclaim{4.5 Corollary} Let $F \to E \to B$ be a fibration in
which $F$ is elliptic with positive Euler characteristic and $B$
is rationally a wedge of odd spheres. If the fibration is TNCZ,
then $E$ is formal, the invariants ${\text{cup}}_0(E)$,
${\text{e}}_0(E)$, ${\text{cat}}_0(E)$ and ${\text{cl}}_0(E)$ all
agree and their common value, ${\text{nil}}_0(E)$, satisfies
${\text{nil}}_0(E) = 1 + {\text{nil}}_0(F)$.
\endproclaim

\demo{Proof}  This follows from Theorem 2.2. \hfill$\square$
\enddemo

So far, we have collected together some consequences of Halperin's
conjecture.  These consequences can be read as weak versions of
Halperin's conjecture. We go on to establish some of these weak
versions of the conjecture.  We can deal quite well with the case
in which the base $B$ is a wedge of odd-dimensional spheres.

The next result generalizes part of
\cite{F\'e-Ha$_2$,Th.10.4(iv)}.

\proclaim{4.6 Proposition} Let $F \to E \to B$ be any fibration in
which $F$ is a rational Poincar\'e duality space and $B$ is a
wedge of odd-dimensional spheres. Then ${\text{e}}_0(E) \geq 1 +
{\text{e}}_0(F)$.
\endproclaim

\demo{Proof} Suppose $\Lambda W, d_B \to \Lambda W\otimes\Lambda
V, D \to \Lambda V,d$ is the minimal model of the fibration.
Observe that $\Lambda W \otimes\Lambda V, D$ must actually be the
minimal model of $E$, i.e., the differential $D$ must be
decomposable. Indeed, this is the case for any fibration in which
$F$ is a space with finite dimensional rational cohomology and $B$
is a wedge of odd-dimensional spheres, as follows from
\cite{Ha$_2$,Th.1.4(iii)}. Observe further, that since $B$ is a
wedge of spheres, it is both formal and coformal.  Thus its
minimal model $\Lambda W, d_B$ is a bigraded model in the sense
discussed earlier, and the differential is quadratic, i.e., $d_B
\: W \to \Lambda^{2}W$.  We use these observations in the proof.

Suppose that ${\text{e}}_0(F) = n$.  Since $F$ is a Poincar\'e
duality space, the fundamental class of $F$ can be represented by
a cocycle $\alpha \in \Lambda^{\geq n} V$ (cf.{}
\cite{F\'e,Lem.5.6.1}). Suppose that $u \in W$ is a generator of
lowest (odd) degree, so that $d_B(u) = 0$.  Consider the element
$u\alpha \in \Lambda^{\geq n+1}(W\oplus V)$. There is no {\it a
priori\/} reason why $u\alpha$ should be a $D$-cocycle, but we
will show the following:

\smallskip

\noindent{\it Claim.}  There is an element $\eta \in \Lambda^{\geq
n+1}(W\oplus V)$, with $\eta \in (\Lambda W)_{+}\otimes\Lambda V$,
such that $D(u\alpha + \eta) = 0$.

\smallskip

\noindent{\it Proof of Claim.} We argue by induction, using the
second grading of the bigraded model $\Lambda W$.  Let $\{ b_{i,
j} \}_{j \in J_i}$ be a (vector space) basis of
$(\Lambda^{+}W)_i$, for each $i\geq 0$.  It is convenient for our
argument to have this basis be a monomial basis, so that each
basis element has a certain length.  Also, our element $u$ is one
of the basis elements $b_{0, j}$, but we denote it $u$ anyway so
as to distinguish it.

Write $D(\alpha) = \sum_{i\geq 0, j} b_{i,j}\,\Omega_{i, j},$ for
suitable elements $\Omega_{i,j} \in \Lambda V$.  Each term
$b_{i,j}\,\Omega_{i, j} \in \Lambda^{\geq n+1}(W\oplus V)$, as $D$
is decomposable.  Then we have
$$D(u\alpha) = -\sum_j u\,b_{0,j}\,\Omega_{0, j}\; -\; \sum_{i\geq
1, j} u\,b_{i,j}\,\Omega_{i, j}.$$
Now $\Lambda W, d_B$ is the minimal model of a wedge of spheres,
whose cohomology has trivial products.  It follows that each
$ub_{0, j} \in (\Lambda ^{\geq 2}W)_0$ is a $d_B$-cocycle.  So
$d_B(\eta_{1, j}) = ub_{0, j}$ for some $\eta_{1, j} \in (\Lambda
W)_1$.  Furthermore, since $\Lambda W, d_B$ is coformal, each
$\eta_{1, j}$ is of length one less than $ub_{0, j}$. For each
$j$, we have
$$D(\eta_{1, j}\,\Omega_{0, j}) = ub_{0, j}\,\Omega_{0,j}\; +\;
(-1)^{|\eta_{1, j}|}\,\eta_{1, j}\,  D(\Omega_{0, j}).$$
From the observation about the length of each $\eta_{1, j}$,
together with the fact that $D$ is decomposable, it follows that
each $\eta_{1, j}\,\Omega_{0, j} \in \Lambda^{\geq n+1}(W\oplus
V)$ and each $\eta_{1, j}\,D(\Omega_{0, j})\in \Lambda^{\geq
n+2}(W\oplus V)$. Finally, write $\eta_{(1)} =  \sum_j \eta_{1,
j}\,\Omega_{0, j}$, so that $\eta_{(1)} \in \Lambda^{\geq
n+1}(W\oplus V)$ and $\eta_{(1)} \in (\Lambda W)_{+}\otimes\Lambda
V$. Then we have
$$D(u\alpha + \eta_{(1)}) = \; -\;\sum_{i\geq 1, j}
u\,b_{i,j}\,\Omega_{i, j} \;+\;\sum_j (-1)^{|\eta_{1,
j}|}\,\eta_{1, j}\,D(\Omega_{0, j})
=
\sum_{i\geq 1, j} b_{i,j}\,\Omega^{(1)}_{i, j},$$
for suitable $\Omega^{(1)}_{i, j} \in \Lambda V$, and each
$b_{i,j} \Omega^{(1)}_{i, j} \in \Lambda^{\geq n+2}(W\oplus V)$.
This starts the induction.

Now suppose inductively that we have an element $\eta_{(r)} \in
\Lambda^{\geq n+1}(W\oplus V)$ with $\eta_{(r)} \in (\Lambda
W)_{+}\otimes\Lambda V$ and $D(u\alpha + \eta_{(r)}) = \sum_{i\geq
r, j} b_{i,j}\, \Omega^{(r)}_{i, j} \in \Lambda^{\geq n+2}(W\oplus
V)$. Re-write $\sum_{j} b_{r,j} \Omega^{(r)}_{r, j}$, the part of
$D(u\alpha + \eta_{(r)})$ that contains terms of lowest second
degree in $\Lambda W$, as follows:  Let $\{c_k\}_{k\in K}$ be a
(vector space) basis for $\Lambda V$.  Then write
$$\sum_{j} b_{r,j}\,\Omega^{(r)}_{r, j}\; =\; \sum_{k}\,\beta_{r,
k}\, c_{k}$$
for suitable terms $\beta_{r, k} \in (\Lambda W)_r$.  For $i\geq
r+1$, we have $D(b_{i,j}\,\Omega^{(r)}_{i, j}) \in (\Lambda
W)_{\geq r}\otimes\Lambda V$. So working modulo the ideal
generated by $(\Lambda W)_{\geq r}$ in $\Lambda W\otimes \Lambda
V$, we have
$$0 = D^2(u\alpha + \eta_{(r)})\; \equiv\;
 \sum_{k}\,d_B(\beta_{r,k})\,c_{k}.$$
Hence $d_B(\beta_{r,k}) = 0$ for each $k$.  Recall once again that
$\Lambda W,d_B$ is the bigraded model, with $H_{+}(\Lambda W,d_B)
= 0$. Since $r\geq 1$, it follows that each $\beta_{r,k}$ is
$d_B$-exact. So $d_B(\eta_{r+1,k}) = \beta_{r,k}$ for some
$\eta_{r+1,k} \in (\Lambda W)_{r+1}$. As $\Lambda W$ is coformal,
$\eta_{r+1,k}c_{k} \in \Lambda^{\geq n+1}(W\oplus V)$ for each
$k$. Now set $\eta_{(r+1)} = \eta_{(r)} - \sum_{k}
\eta_{r+1,k}c_{k}$. Note that $\eta_{(r+1)} \in \Lambda^{\geq
n+1}(W\oplus V)$ and $\eta_{(r+1)} \in (\Lambda
W)_{+}\otimes\Lambda V$. A straightforward check shows that
$$\align D(u\alpha + \eta_{(r+1)})\; &=\;
 \sum_{i\geq r+1, j} b_{i,j}\,\Omega^{(r)}_{i, j}
\;-\;\sum_{k}\,(-1)^{|\eta_{r+1, k}|}\,\eta_{r+1,k}\,D(c_{k}) \\
&= \sum_{i\geq r+1, j}\,b_{i,j}\,\Omega^{(r+1)}_{i, j},\\
\endalign
$$
for suitable terms $\Omega^{(r+1)}_{i, j} \in \Lambda V$ with each
$b_{i,j}\,\Omega^{(r+1)}_{i, j} \in \Lambda^{\geq n+2}(W\oplus
V)$. This completes the inductive step.

Since $B$ is simply connected, the lowest degree of a generator in
$W_r$ increases strictly with $r$.  Hence, by taking $r$
sufficiently large, we obtain an element $\eta_{(r)}$ as in the
claim, with $D(u\alpha + \eta_{(r)}) = 0$. {\it End of Proof of
Claim.\/}

We now show that this cocycle is not $D$-exact.  Suppose that
$D(\zeta + \chi) = u\alpha +\eta$ for $\zeta \in \Lambda V$ and
$\chi \in \Lambda^{+} W\otimes\Lambda V$. Then $d(\zeta) = 0$.
However, $\zeta$ has higher degree than $\alpha$, which represents
the fundamental class of $F$. Thus $\zeta = d(a)$ for some $a \in
\Lambda V$ so $D(a) = \zeta + \chi'$ for some $\chi' \in
\Lambda^{+} W\otimes\Lambda V$. Write $\chi - \chi' = \sum_{i,
j}b_{i,j}\chi_{i,j}$.  Working modulo the ideal in $\Lambda
W\otimes\Lambda V$ generated by $\Lambda^{\geq2}W_0 + (\Lambda
W)_{+}$, we have $D(\chi+\zeta) = D(\chi - \chi') =
\sum_{j}(-1)^{|b_{0,j}|}\,b_{0,j}\,d(\chi_{0,j}) \equiv u\alpha$.
This implies $\alpha$ is $d$-exact, which is a contradiction since
$\alpha$ represents the fundamental class of $F$.  Therefore,
$u\alpha + \eta$ is a non-exact $D$-cocycle in $\Lambda^{\geq
n+1}(W\oplus V)$.  The result follows. \hfill$\square$
\enddemo

Next we give the main result of this section.

\proclaim{4.7 Theorem} Let $F \to E \to B$ be a fibration in which
$F$ is elliptic with positive Euler characteristic and $B$ is a
wedge of odd-dimensional spheres.  Then for $E$ we have
${\text{e}}_0(E) = {\text{cat}}_0(E) = {\text{cl}}_0(E)$ and
furthermore these equal ${\text{nil}}_0(F)+1$.
\endproclaim

\demo{Proof} In fact the proof will display a simple model for $E$
of homotopical nilpotency equal to $1 + {\text{nil}}_0(F)$.  Let
$\rho: \Lambda V, d \to H(F)$ be the bigraded model for $F$. By
construction we have $\rho(V_1) = 0$ and hence $\rho\big((\Lambda
V)_+\big) = 0$.  Consider the ideal ${\text{ker}}\,\rho
\subseteq\Lambda V$. This is a differential ideal since all
boundaries are in ${\text{ker}}\,\rho$. Further, it is an acyclic
ideal since $\rho$ is a surjective quasi-isomorphism. Since $B$ is
formal, the fibration has a model $H(B) \to H(B)\otimes\Lambda V,
D \to \Lambda V, d$. We claim that the ideal
$H(B)\otimes{\text{ker}}\,\rho$ of $H(B)\otimes\Lambda V$ is also
a differential acyclic ideal.  To see this, we argue as follows:

Let $\{b_i\}$ be a  basis for $H^+(B)$. For any $\chi \in \Lambda
V$,  we can write
$$D(\chi) = d\chi + \sum_i\,b_i \Omega_i(\chi),$$
and a standard argument shows that this defines derivations
$\Omega_i$ on $\Lambda V, d$, each of negative even degree
according as the degree of the $b_i$ (cf. the result of Meier,
cited as Theorem 1.5). Now for $w \in V_1$, we have $D(w) = dw +
\sum_i\,b_i \Omega_i(w)$. For parity of degree reasons we must
have $\Omega_i(w) \in (\Lambda V)_{+}$.  It follows that
$\Omega_i\big((\Lambda V)_+\big) \subseteq (\Lambda V)_+$ for each
derivation $\Omega_i$.

Next, if $\chi \in {\text{ker}}\,\rho\subseteq \Lambda V$, write
$\chi = \chi_0 + \chi_{+}$ where $\chi_0 \in \Lambda V_0$ and
$\chi_{+} \in (\Lambda V)_+$. Since $\chi \in {\text{ker}}\,\rho$,
it follows that $\chi_0 \in {\text{ker}}\,\rho$ and hence $\chi_0
= d\eta$ for some $\eta \in (\Lambda V)_1$.  Then
$$ \align D(\chi) &= d(\chi_{+}) + \sum_i\,b_i \Omega_i(\chi_0) +
\sum_i\,b_i \Omega_i(\chi_{+}) \\ &= d(\chi_{+}) + \sum_i\,b_i
\Omega_i(d\eta) + \sum_i\,b_i \Omega_i(\chi_{+}) \\ &= d(\chi_{+})
+ \sum_i\, b_i d\big(\Omega_i(\eta)\big) + \sum_i\,b_i
\Omega_i(\chi_{+}). \\
\endalign
$$
This term is in $H(B)\otimes{\text{ker}}\,\rho$ because
${\text{ker}}\,\rho$ contains all boundaries in $\Lambda V$ and
also, as remarked earlier, $\Omega_i(\chi_{+}) \in
\Omega_i\big((\Lambda V)_+\big) \subseteq (\Lambda V)_+ \subseteq
{\text{ker}}\,\rho$.   We have shown that $D({\text{ker}}\,\rho)
\subseteq H(B)\otimes{\text{ker}}\,\rho$ and since $D=0$ on
$H(B)$, it follows that $H(B)\otimes{\text{ker}}\,\rho$ is
$D$-stable.

Further, suppose $\alpha \in H(B)\otimes{\text{ker}}\,\rho$ is a
$D$-cocycle. We can write $\alpha = \alpha' + \sum_i\,b_i
\alpha_i$ with $\alpha'$ and each $\alpha_i$ in
${\text{ker}}\,\rho$. Then $D(\alpha) = 0$ implies $d(\alpha') =
0$ and hence $\alpha' = d\eta'$ for some $\eta' \in
{\text{ker}}\,\rho$, as ${\text{ker}}\,\rho$ is acyclic. Without
loss of generality, we can choose $\eta' \in (\Lambda V)_+$, as
$d(V_0)=0$. Then $D(\eta') = d(\eta') +
\sum_i\,b_i\Omega_i(\eta')$. From an earlier remark, each
$\Omega_i(\eta') \in (\Lambda V)_+ \subseteq {\text{ker}}\,\rho$.
So we have
$$\alpha = D(\eta') + \sum_i\,b_i\big(\alpha_i -
\Omega_i(\eta')\big) = D(\eta') + \sum_i\,b_i \alpha'_i$$
for elements $\eta', \alpha'_i \in {\text{ker}}\,\rho$. As all
products in $H^+(B)$ are trivial, $D(\alpha) = 0$ implies each
$d(\alpha'_i) = 0$, so that $\alpha'_i = d\eta_i$ for elements
$\eta_i \in {\text{ker}}\,\rho$. But $D(b_i \eta_i) = - b_i
d(\eta_i)$ and hence $\alpha = D(\eta' - \sum_i\,b_i \eta_i)$ with
$\eta' - \sum_i\,b_i \eta_i \in H(B)\otimes{\text{ker}}\,\rho$.
This shows  $H(B)\otimes{\text{ker}}\,\rho$ is an acyclic ideal of
$H(B)\otimes\Lambda V$.

To finish, notice that the projection
$$q: H(B)\otimes\Lambda V \to {H(B)\otimes\Lambda V \over
H(B)\otimes{\text{ker}}\,\rho}$$
is a quasi-isomorphism, as $H(B)\otimes{\text{ker}}\,\rho$ is
acyclic. If ${\text{nil}}_0(F) = n$, then $\Lambda^{>n} V
\subseteq \Lambda^{>n} V_0 + (\Lambda V)_+ \subseteq
{\text{ker}}\,\rho$ and hence the quotient ${H(B)\otimes\Lambda V
\slash(H(B)\otimes{\text{ker}}\,\rho})$ is a nilpotent DG algebra
of length $\leq n+1$, and is quasi-isomorphic to a KS-model for
$E$. Therefore, it follows from \cite{Co} that ${\text{cl}}_0(E)
\leq n + 1$. From Proposition 4.6, we have $n+1 \leq
{\text{e}}_0(E)$.  But in general we have ${\text{e}}_0(E) \leq
{\text{cat}}_0(E) \leq {\text{cl}}_0(E)$.  Hence these three
invariants must agree and furthermore must equal $n+1$.
\hfill$\square$
\enddemo

\heading \S5 --- Discussion and Questions
\endheading

In this last, somewhat informal, section we discuss the foregoing
results and give some additional ones.  The intention is to
indicate both limits on, and possibilities for, development of the
work in Sections 3 and 4. We also offer some specific questions
along these lines. All this is collected under two subheadings,
according as the topic relates to Section 3 or Section 4.

\subheading{5.1 Formality}

One can attempt a quite general study of relations between the
formality of $F$, $E$ and $B$ in a fibration (cf. \cite{Th$_2$}
and \cite{Vi}). It is tempting to think that the very restrictive
hypotheses of Proposition 3.2 might be relaxed, whilst keeping the
conclusion.  However, the following sort of example suggests these
hypotheses are actually sharp.

\smallskip

\noindent{\bf Example.}  There is a rational fibration $S^2 \vee
S^2 \vee S^2 \to E \to S^3$ that is TNCZ in which $E$ is not
formal.  We omit details of this example, as it is similar to
\cite{Th$_2$,Ex.III.13} and other examples. We note that the fibre
is formal with non-zero cohomology only in even degrees. The
fibration we have in mind is actually much closer to being trivial
than just TNCZ.  It satisfies $H^*(E;{\Bbb Q}) \cong H^*(B;{\Bbb
Q})\otimes H^*(F;{\Bbb Q})$ as algebras and admits a rational
section.

In case a weaker conclusion is acceptable, then of course there
are more possibilities.  In Proposition 3.2, we can remove the
hypothesis of ellipticity on the fibre space, but so far only at
the cost of a greatly weakened conclusion, thus:

\proclaim{Proposition}  Let $F \to E \to B$ be a fibration in
which $F$ and $B$ are both formal. If the fibration is TNCZ, then
$E$ has spherically-generated cohomology.
\endproclaim

\noindent{}The proof of this proposition is omitted. It can be
proved with an argument similar to that of Proposition 3.2.  The
above example also suggests this might be a sharp conclusion,
without much more restrictive hypotheses.

It would be satisfying to have a converse of Proposition 3.2.  In
Proposition 3.3, we have such a result but only for a very special
base.  We offer the following as a specific question in this area:

\proclaim{Question 1}  Let $F \to E \to B$ be a fibration in which
$F$ is  formal and elliptic and $B$ is formal. If $E$ is formal,
then is the fibration is TNCZ?
\endproclaim

Of course, there are many variations on this type of question to
investigate. For more results and examples along these lines, see
\cite{Th$_2$} and \cite{Vi}.

Returning to Conjecture 1.1, it may be possible to apply Theorem
3.4, together with an appropriate obstruction theory for the
formality of $E$ in such a fibration. Such an obstruction theory
does exist, and I hope to develop this line of investigation in
future work.

\subheading{5.2 Numerical Invariants}

The results of Section 4 deal nicely with fibrations $F \to E \to
B$, in which $F$ is positively elliptic and $B$ is a wedge of
odd-dimensional spheres, but only for the invariants
${\text{e}}_0$, ${\text{cat}}_0$ and ${\text{cl}}_0$.

\proclaim{Question 2} For a fibration $F \to E \to B$, in which
$F$ is elliptic with positive Euler characteristic and $B$ is a
wedge of odd-dimensional spheres, is ${\text{cup}}_0(E) = 1 +
{\text{nil}}_0(F)$?
\endproclaim

Of course, if Halperin's conjecture is true, it implies an
affirmative answer to Question 2 (Corollary 4.5).   It is curious
that Question 2 remains unresolved whilst its counterpart for the
invariants ${\text{e}}_0$, ${\text{cat}}_0$ and ${\text{cl}}_0$ is
established.

Another `asymmetry' in the results we have presented is the lack
of an inequality for cone-length in Proposition 4.3.  This gives
our next specific question:

\proclaim{Question 3} Let $F \to E \to B$ be a fibration with $F$
formal.  If the fibration is TNCZ, then is ${\text{cl}}_0(E) \geq
{\text{cl}}_0(B) + {\text{nil}}_0(F)$?
\endproclaim

As regards extending the results presented in Section 4, we can
show at least the following:

\proclaim{Proposition} Let $F \to E \to B$ be a fibration in which
$F$ is elliptic with positive Euler characteristic and $B$ has
${\text{e}}_0(B) = 1$. Then ${\text{e}}_0(E) \geq 1 +
{\text{nil}}_0(F)$.
\endproclaim

Of course, if ${\text{e}}_0(B) = 1$, then the base space $B$ is
rationally a wedge of spheres and ${\text{e}}_0(B) =
{\text{cat}}_0(B) = {\text{cl}}_0(B) = 1$.  Notice, however, that
we allow even dimensional spheres and also that we do not have
Poincar\'e duality in either the base or the total spaces.  The
current proof of this Proposition uses a rather involved reduction
argument, similar to that of \cite{Th$_1$} or Proposition 4.6
here. This result suggests that the `next step' might be to
investigate one of the following questions:

\proclaim{Question 4 (a)} Let  $F \to E \to B$ be a fibration in
which $F$ is elliptic with positive Euler characteristic and
${\text{e}}_0(B) = 2$.  Is ${\text{e}}_0(E) \geq 2 +
{\text{nil}}_0(F)$?
\endproclaim

\proclaim{Question 4 (b)} Let  $F \to E \to B$ be a fibration in
which $F$ is elliptic with positive Euler characteristic.  Is
${\text{e}}_0(E) \geq 1 + {\text{nil}}_0(F)$?
\endproclaim

In the previous two paragraphs, we focussed on the
${\text{e}}_0$-invariant.  Choosing one of the other three
invariants gives similar questions to be investigated. There is
some overlap between all these questions, because of the results
of \cite{Co-F\'e-Le} and \cite{F\'e-Ha-Le} mentioned earlier:
${\text{e}}_0(X) = {\text{cat}}_0(X) = {\text{cl}}_0(X)$ whenever
$X$ is a rational Poincar\'e duality space. Thus, if $F \to E \to
B$ is a fibration in which $F$ and $B$ are rational Poincar\'e
duality spaces, then so too is $E$ a rational Poincar\'e duality
space and hence the three invariants agree on each of $F$, $E$ and
$B$.  Notice that this observation obtains the first part of the
conclusion to Theorem 4.7 in the case $B = S^{2n+1}$.  Finally, we
mention that Jessup's result \cite{Je,Prop.3.6} also shows that
${\text{cat}}_0(E) \geq {\text{cat}}_0(B) +1$ for any fibration in
which $F$ is positively elliptic --- without the hypothesis of
TNCZ.  These latter observations give some interesting complements
to the results presented here.

\enddocument